\documentclass{article}
\usepackage{amsmath}
\usepackage{amsfonts}
\usepackage{amssymb}
\usepackage{amsthm}
\usepackage{amscd}
\usepackage{makeidx}
\usepackage{a4wide}
\makeindex
\newtheorem{theorem}{Theorem}
\newtheorem{remark}{Remark}
\newtheorem{proposition}{Proposition}
\newtheorem{definition}{Definition}
\newtheorem{lemma}{Lemma}
\newtheorem{cor}{Corollary}

\newcommand{\supp}{\text{supp }}
\newcommand{\e}{\varepsilon}
\begin{document}
\date{}

\title{Beurling-Pollard type theorems}
\author{{\it Victor Shulman}  and
{\it Lyudmila  Turowska} }


\maketitle{}
\begin{abstract}
We establish a version of the Beurling-Pollard theorem for operator
synthesis and apply it to derive some results on linear operator
equations and to prove a Beurling-Pollard type theorem for Varopoulos
tensor algebras. Additionally we establish a Beurling-Pollard theorem
for weighted Fourier algebras and use it to obtain ascent estimates for operators that are functions of generalized scalar operators.
\end{abstract}

 \section{Introduction.}

The Beurling - Pollard Theorem is the result (or a group of results,
see
\cite{beurling, pollard, kahane, kinukawa, walker})
stating that a function in the Fourier algebra of a locally compact
Abelian group admits the spectral synthesis if it is sufficiently
"smooth", or (and) its null set is sufficiently thin. Varopoulos
\cite{var} obtained a result of this kind for the tensor algebras
$C(X)\hat{\otimes}C(Y)$.

Recall that by spectral synthesis in a commutative semisimple regular Banach algebra (or more generally, in a normed
space of functions) one understands a possibility to approximate functions equal to zero on some set by functions
equal to zero on a neighborhood of this set.

Our first aim is to obtain a version of Beurling - Pollard theorem
for operator synthesis. The latter means a synthesis in the special
space of functions on the product of two measure spaces: the space
of integral kernels of nuclear operators between $L^2$-spaces. It
should be noted that the measure spaces are not supplied with
topologies, but their product carries a natural {\it
pseudo-topology} and the "neighborhoods" are taken with respect to
it. The importance of the operator synthesis for the theory of
operator algebras and their invariant subspaces was discovered by
Arveson \cite{arv}.

We obtain the operator synthesis B-P theorem in Section 2. Then in
Section 3 we apply it to the linear operator equations with normal
coefficients. This allows us to obtain easily some results in
\cite{sht2} and extend them to more general situations. More
specifically, the main result of Section 2 gives an estimation of
ascent for general "normal" multiplication operators via the
Hausdorff dimension of one of the coefficient families. Moreover the
estimates of the ascent are given in their strong form: as the
statement on coincidence of the root space of an operator with the
kernel of some of its power.

It should be noted that the estimation of the ascent of linear multiplication operators is
very important for the theory of linear operator equations (see for
example \cite{sht2}). Moreover the classical forms of the spectral
synthesis problems in group algebras can be reformulated in terms of the ascent of
multiplication operators on the dual spaces of the algebras.

In Section 4, using established in \cite{sht1} connection between
operator synthesis and spectral synthesis in Varopoulos algebras, we
obtain a strengthened version of the B-P theorem of Varopoulos
\cite{var} and show that from this form of Varopoulos's theorem one
can easily deduce the classical B-P theorem for Fourier algebras.

In Section 5 we prove a B-P theorem for weighted Fourier algebras
$A_{\alpha}({\mathbb T}^n)$ on $n$-dimensional tori. In this theorem
the conditions that provide synthesis for a function $f\in
A_{\alpha}({\mathbb T}^n)$, depend on {\it balanced Hausdorff
measure} (the notion introduced in \cite{keckic}) of the null set of
$f$.

In the last section we apply the results of Section 5 to obtain
estimates of the ascent for operators that are expressed as
functions of commuting generalized scalar operators. As a special
case this gives us the bounds for ascent of an elementary operators
$x\to \sum_{i=1}^na_ixb_i$  obtained previously in \cite{keckic}.

We are indebted to Milos Arsenovic and Dragolyub Keckic for sending
us a preprint of their very interesting paper \cite{keckic}.

\section{The Beurling-Pollard Theorem for operator synthesis.}

First we recall some definitions and facts from
\cite{arv,sht1,sht2}. Let $(X,\mu)$, $(Y,\nu)$ be standard measure
spaces with $\sigma$-finite measures. We say that a measurable
subset $M$ of $X\times Y$ is {\it marginally null} if $M\subset
(X_1\times Y)\cup(X\times Y_1)$, where $\mu(X_1)=\nu(Y_1)=0$. If the
symmetric difference of subsets $M$, $N$ is marginally null we say
that $M$, $N$ are marginally equivalent and write $M\simeq N$.
Following \cite{eks} we define an $\omega$-topology on $X\times Y$
such that the $\omega$-open (pseudo-open) sets are, modulo
marginally null sets, countable unions of measurable rectangles,
i.e. the sets $\alpha\times\beta$ with measurable $\alpha$, $\beta$.
The complements to $\omega$-open sets are called $\omega$-closed or
pseudo-closed sets.

Let $\Gamma(X,Y)=L_2(X,\mu)\hat\otimes L_2(Y,\nu)$, the projective
tensor product of $L_2$-spaces. Each $\Psi\in \Gamma(X,Y)$ can be
identified with a function $\Psi:X\times Y\to{\mathbb
C}$ which admits a representation
\begin{equation}\label{rep}
\Psi(x,y)=\sum_{n=1}^{\infty}f_n(x)g_n(y), \end{equation} such that
$f_n\in L_2(X,\mu)$, $g_n\in L_2(Y,\nu)$ and
$\sum_{n=1}^{\infty}||f_n||_{L_2}\cdot ||g_n||_{L_2}<\infty$. Such a
function is defined marginally almost everywhere (that is uniquely up to a marginally null set).
Moreover,
$$||\Psi||_{\Gamma}=\text{inf}\sum_{n=1}^{\infty}||f_n||_{L_2}\cdot
||g_n||_{L_2},$$ where infimum is taken over all representations
(\ref{rep}) of $\Psi$.

By \cite{eks}, each function $\Psi\in \Gamma(X,Y)$ is
pseudo-continuous, i.e. the preimage of any open subset of ${\mathbb
C}$ is pseudo-open. We say that $\Psi$ vanishes on $E\subset X\times
Y$ if $F\chi_E=0$ m.a.e.\, where $\chi_E$ is the characteristic
function of $E$. For ${\mathcal F}\subset \Gamma(X,Y)$, the null
set, $\text{null }{\mathcal F}$, is defined to be the largest, up to
a marginally null set, pseudo-closed set such that each function
$F\in{\mathcal F}$ vanishes on it. For a pseudo-closed set $E$,
set\begin{eqnarray*}
&\Phi(E)=\{\Psi\in \Gamma(X,Y): \Psi\text{ vanishes on }E\},\\
&\Phi_0(E)=\overline{\{\Psi\in \Gamma(X,Y): \Psi\text{ vanishes on a
pseudo-nbhd of } E\}}.
\end{eqnarray*}
Then $\Phi(E)$ and $\Phi_0(E)$ are closed subspaces of
$\Gamma(X,Y)$, invariant under  multiplication by
$L_{\infty}(X,\mu)$ and $L_{\infty}(Y,\nu)$ functions and such that
$\text{null }\Phi(E)=\text{null }\Phi_0(E)=E$. Moreover, if
$A\subset \Gamma(X,Y)$ is another invariant closed subspace such
that $\text{null }A=E$ then $\Phi_0(E)\subset A\subset \Phi(E)$
(\cite[Theorem~2.1]{sht1}).

A subset $E\subset X\times Y$ is called {\it a set of
operator synthesis} (or {\it an operator synthetic set}) (with respect to
$(\mu,\nu)$) if $\Phi(E)=\Phi_0(E)$.

The space $\Gamma(X,Y)$ is  predual to the space of bounded
operators $B(H_1,H_2)$ from $H_1=L_2(X,\mu)$ to $H_2=L_2(Y,\nu)$.
The duality is given by
$$\langle X,\Psi\rangle=\sum_{n=1}^{\infty} (Xf_n,\bar g_n),$$
where $X\in B(H_1,H_2)$ and $\Psi(x,y)=\sum_{n=1}^{\infty} f_n(x)
g_n(y)$. It can be said that $\Gamma(X,Y)$ is the space of integral kernels of nuclear operators.

Let $P(U)$ and $Q(V)$ denote the multiplication operators by the
characteristic functions of $U\subset X$ and $V\subset Y$. We say
that $T\in B(H_1,H_2)$ is {\it supported} in $E\subset X\times Y$
(or $E$ supports $T$) if $Q(V)TP(U)=0$ for each measurable sets
$U\subset X$, $V\subset Y$ such that $(U\times V)\cap E\simeq
\emptyset$.
 Then there exists the smallest (up to a marginally null set)
pseudo-closed set, $\supp T$, which supports $T$. More generally,
for any subset ${\mathfrak M}\subset B(H_1,H_2)$ there is the
smallest pseudo-closed set, $\supp{\mathfrak M}$, which supports all
operators in ${\mathfrak M}$ (see \cite{sht1}).

For any pseudo-closed set $E\subset X\times Y$, the set ${\mathfrak
M}_{max}(E)$ of all operators $T$, supported in $E$, has support
$E$.  It is easy to check that ${\mathfrak M}_{max}(E)$ is a
${\mathcal D}_1\times{\mathcal D}_2$-bimodule, where ${\mathcal
D}_1$, ${\mathcal D}_2$  are the algebras of multiplications by
functions in  $L_{\infty}(X,\mu)$ and $L_{\infty}(Y,\nu)$
respectively. Clearly ${\mathfrak M}_{max}(E)$ is the largest
$\sigma$-closed bimodule with support equal to $E$. There is also
the smallest $\sigma$-closed bimodule ${\mathfrak M}_{min}(E)\subset
B(H_1,H_2)$ with support equal to  $E$ (\cite{arv}), and moreover
$${\mathfrak M}_{max}(E)=\Phi_0(E)^{\perp}, \quad
{\mathfrak M}_{min}(E)=\Phi(E)^{\perp}$$ (see \cite{sht1}). So
a pseudo-closed set $E\subset X\times Y$ is  operator synthetic
 iff  one of the following equivalent conditions holds:
\begin{itemize}
\item ${\mathfrak M}_{max}(E)={\mathfrak M}_{min}(E)$.
\item $\langle X,\Psi\rangle=0$ for any $X\in B(H_1,H_2)$ and
$\Psi\in\Gamma(X,Y)$ with $\supp X\subset E\subset\text{null }\Psi$.
\end{itemize}

We will also consider the space $V_{\infty}(X,Y)=L^{\infty}(X,\mu)\hat\otimes
L^{\infty}(Y,\nu)$ of all (marginal
equivalence classes of) functions $\Psi(x,y)$ that can be written in
the form (\ref{rep}) with $f_n\in L^{\infty}(X,\mu)$, $g_n\in
L^{\infty}(Y,\nu)$ and
$$\sum_{n=1}^{\infty} ||f_n||_{\infty}^2\leq C,\quad
\sum_{n=1}^{\infty}||g_n||_{\infty}^2\leq C.$$ The least possible
$C$ here is the norm of $\Psi$ in $V_{\infty}(X,Y)$. We have that
$V_{\infty}(X,Y)\Gamma(X,Y)\subset\Gamma(X,Y)$.

The space $B(H_1,H_2)$ is a $V_{\infty}(X,Y)$-module: for
$F=\sum_{n=1}^{\infty} f_n\hat\otimes g_n\in V^{\infty}(X,Y)$ and
$X\in B(H_1,H_2)$, we set
$$F\cdot T=\sum_{n=1}^{\infty}M_{g_n}XM_{f_n},$$
where $M_f$ is the multiplication operator by $f$. Here the sum
converges in the norm operator topology.

We denote by $\Delta_F$
the operator $X\mapsto F\cdot X$ on $B(H_1,H_2)$. Let us say that $F\in
V_{\infty}(X,Y)$ {\it  admits operator synthesis} (is {\it operator
synthetic}) with respect to $(\mu, \nu)$ if
$$\ker\Delta_F={\mathfrak M}_{max}(\text{null }F).$$
Observe that the inclusion $\ker\Delta_F\subset{\mathfrak
  M}_{max}(\text{null }F)$ holds for any $F\in V_{\infty}(X,Y)$, see
\cite[Lemma~4.1]{sht2}.

Note that if measures $\mu$, $\nu$ are finite then $
V_{\infty}(X,Y)\subset \Gamma(X,Y)$. In this case the notion agree
with the one in spectral synthesis: $F\in V_{\infty}(X,Y)$ admits
operator synthesis with respect to ($\mu$, $\nu$) iff $F\in
\Phi_0(\text{null }F)$, i.e.
$$\langle T,F\rangle=0\text{ whenever }\supp T\subset\text{null }F.$$
Indeed the equality $\langle T, F\rangle=\langle F\cdot T, 1\rangle$  shows that any operator synthetic element $F$ of $ V_{\infty}(X,Y)$ belongs to $\Phi_0(\text{null }F)$. Conversely, if $F\in\Phi_0(\text{null
}F)\cap V_{\infty}(X,Y)$
then $F\Psi\in \Gamma(X,Y)\cap \Phi_0(\text{null F})$ for any $\Psi\in\Gamma(X,Y)$, and
$\langle F\cdot T,\Psi\rangle=\langle T, F\Psi\rangle=0$ whenever
$\supp T\subset\text{null }F$, implying $F\cdot T=0$.

From now on we will assume  that $X$, $Y$ are compacts, $X$ is a
metric space with
 metric $d$ and the Hausdorff dimension $w<\infty$, and $\mu$, $\nu$
 are $\sigma$-finite regular Borel measures on $X$ and $Y$ respectively.
For $x\in X$ and $A\subset X$ we  denote by $d(x,A)$ the distance
between $A$ and $x$, i.e. $d(x,A)=\text{inf}_{y\in A} d(x,y)$; we
assume $d(x,\emptyset)=\infty$.

\begin{lemma} \label{lemma1} Let $E\subset X\times Y$ such that $E^c$ is a countable union of
Borel rectangles. Then $E$ is a union of countable number of compact
sets and a marginally null set.
\end{lemma}

\begin{proof}
Let $\displaystyle E^c=\bigcup_{n=1}^{\infty} A_n\times B_n$. Given
$\varepsilon>0$ there exist open sets $A_n^{\varepsilon}\supset A_n$
and $B_n^{\varepsilon}\supset B_n$ such that
$\mu(A_n^{\varepsilon}\setminus A_n)<\e/2^n$ and
$\nu(B_n^{\e}\setminus B_n)<\e/2^n$ and thus $(A_n^{\e}\times
B_n^{\e})\setminus (A_n\times B_n) \subset X_{\e}^n\times Y\cup
X\times Y_{\e}^n$ with $\mu(X_{\e}^n)<\e/2^n$ and
$\nu(Y_{\e}^n)<\e/2^n$ ($X_{\e}^n=A_n^{\varepsilon}\setminus A_n$,
$Y_{\e}^n=B_n^{\e}\setminus B_n$).

Set $E_{\e}=\displaystyle (\bigcup_{n=1}^{\infty} A_n^{\e}\times
B_n^{\e})^c$. Then $E_{\e}$ is compact, $E_{\e}\subset E$ and
 $E\setminus E_{\e}\subset (X_{\e}\times
Y)\cup (X\times Y_{\e})$, where $\displaystyle
X_{\e}=\bigcup_{n}X_{\e}^n$, $Y_{\e}=\bigcup_n Y_{\e}^n$. Take a
decreasing to zero sequence $\{\e_k\}$.  Choosing $A_n^{\e_k}$,
$B_n^{\e_k}$ so that $A_n^{\e_k}\supset A_n^{\e_{k+1}}$,
$B_n^{\e_k}\supset B_n^{\e_{k+1}}$  we obtain $X_{\e_k}\supset
X_{\e_{k+1}}$, $Y_{\e_k}\supset Y_{\e_{k+1}}$ and
$$E\setminus
(\bigcup_k E_{\e_k})\subset (\bigcap_k X_{\e_k}\times
Y)\bigcup(X\times \bigcap_k Y_{\e_k}).$$ As $\mu (X_{\e_k})<\e_k$
and $\nu(Y_{\e_k})<\e_k$, $\displaystyle E\setminus (\bigcup_k
E_{\e_k})$ is marginally null.
\end{proof}

For a subset $E$ of  $X\times Y$ and $y\in Y$, let
$$E^y=\{x\in X: (x,y)\in E\},$$
the {\it $X$-section} of $E$, defined by $y$.
Then clearly, $d(x,E^y)=\text{inf}\{d(x,x_1): (x_1,y)\in E\}$.

\begin{theorem}\label{opsyn}
Let $E$ be a pseudo-closed set and let $F\in V_{\infty}(X,Y)$  satisfy the condition
\begin{equation}\label{eq1}
|F(x,y)|\leq Cd(x,E^y)^{\rho}, \quad \text{for some }C>0.
\end{equation}
Assume that $\rho\geq w/2$. Then ${\mathfrak M}_{max}(E)\subset
\ker\Delta_F$. In particular, if $E=\text{null } F$ then $F$ obeys
operator synthesis.
\end{theorem}

\begin{proof} Our first step is to prove that $\Delta_{F}(T)$ is a Hilbert-Schmidt operator for
any $T\in {\mathfrak M}_{max}(E)$.

For $\varepsilon >0$ let ${\mathcal
E}=(\alpha_1,\ldots,\alpha_{n(\varepsilon)})$ be a family of
pairwise disjoint closed subsets of $X$ such that  $\text{diam
}\alpha_k\leq\varepsilon$ and $\mu (X\setminus \cup\alpha_i)<\e$,
and let $e_k=\chi_{\alpha_k}/||\chi_{\alpha_k}||$, where
$\chi_{\alpha_k}$, the characteristic function of the set
$\alpha_k$.

Denote by
$P_{\mathcal
  E}$ the projection onto the (finite-dimensional) subspace spanned by all $e_k$. Then
$TP_{\mathcal E}\in {\mathfrak S}_2$, $TP_{\mathcal E}\to^s T$, as
$\varepsilon\to 0$ and $\displaystyle\langle F\cdot
T,\Psi\rangle=\lim_{\varepsilon\to 0}\langle TP_{\mathcal
E},F\Psi\rangle$ for each $\Psi\in \Gamma(X,Y)$.

Note that the condition (2) will be preserved if we delete from $X$
and $Y$ some subsets.
Moreover, if for a given $\delta>0$, the sets $X_{\delta}\subset X$, $Y_{\delta}\subset Y$ are such
that $\mu(X_{\delta})<\delta$, $\nu(Y_{\delta})<\delta$ then for
$T_{\delta}:=Q(Y\setminus Y_{\delta})TP(X\setminus Y_{\delta})$ we
have $\supp T_{\delta}\subset E\cap (X\setminus
X_{\delta}\times(Y\setminus Y_{\delta})$, $T_{\delta}\to T$
strongly as $\delta\to 0$. Therefore if we prove that
$\Delta_F(T_{\delta})=0$ for
each $\delta>0$, then it will imply that $\Delta_F(T)=0$.
Taking this into account we may delete from $X$
and $Y$ some open  subsets of small measures and assume (since $E$ is a pseudo-closed
set) that $E^c$ is a countable union of Borel rectangles. Moreover,
by Lemma~\ref{lemma1}, we may repeat this trick and suppose that
$E=\cup_{n=1}^{\infty}E_n$ for some compact sets  $E_n$.

For a
closed subset $\alpha\subset X$,
let $\beta(\alpha)=\pi_2(E\cap(\alpha\times Y))$, where
$\pi_2$ is the projection on the second coordinate. Then
$\beta(\alpha)$ is Borel as a countable union of compact sets.
 As $T$ is supported in $E$, $TP(\alpha)=Q(\beta(\alpha))TP(\alpha)$. Thus
\begin{eqnarray}
|\langle TP_{\mathcal E}, F\Psi\rangle|
&=&|\langle\sum_{i=1}^{n(\varepsilon)}Q(\beta(\alpha_i))
TP_{\mathcal
E}P(\alpha_i),F\Psi\rangle|\nonumber\\&=&|\sum_{i=1}^{n(\varepsilon)}\langle
TP_{\mathcal
E}P(\alpha_i),\chi_{\alpha_i}(s)\chi_{\beta(\alpha_i)}(t)F(s,t)\Psi(s,t)\rangle|
\nonumber\\&\leq& \sum_{i=1}^{n(\varepsilon)}||TP_{\mathcal
E}P(\alpha_i)||_{{\mathfrak S}_2}||\chi_{\alpha_i}(s)
\chi_{\beta(\alpha_i)}(t)F(s,t)\Psi(s,t)||_{L_2}\nonumber\\
&\leq&
\sum_{i=1}^{n(\varepsilon)}||Te_i||\cdot||\chi_{\alpha_i}(s)\chi_{\beta(\alpha_i)}(t)
F(s,t)||_{L_{\infty}} ||\chi_{\alpha_i}(s)\Psi(s,t)||_{L_2}\nonumber\\
&\leq&
||T||\left(\sum_{i=1}^{n(\varepsilon)}||\chi_{\alpha_i}(s)\chi_{\beta(\alpha_i)}(t)
F(s,t)||_{L_{\infty}}^2\right)^{1/2}\left(\sum_{i=1}^{n(\varepsilon)}
||\chi_{\alpha_i}(s)\Psi(s,t)||_{L_2}^2\right)^{1/2}\nonumber\\
&\leq
&||T||\left(\sum_{i=1}^{n(\varepsilon)}||\chi_{\alpha_i}(s)\chi_{\beta(\alpha_i)}(t)
F(s,t)||_{L_{\infty}}^2\right)^{1/2}||\Psi||_{L_2}\label{last}
\end{eqnarray}
Choosing, for each $t\in\beta(\alpha_i)$, a point $s(t)\in
\alpha_i$, one has
\begin{eqnarray}\label{F(s,t)}
\displaystyle |\chi_{\alpha_i}(s)\chi_{\beta(\alpha_i)}(t)F(s,t)|
\leq C d(s,E^t)^{\rho}\leq C d(s,s(t))^{\rho}\leq C(\text{diam
}\alpha_i)^{\rho} \quad\text {a.e.} \end{eqnarray} As the Hausdorff
dimension of $X$ is $w$, we can choose $\alpha_i$ in such a way that
$\sum_{i=1}^{n(\e)}(\text{diam }\alpha_i)^{w}\leq K$, for some
constant $K$ which does not depend on $\e$. Moreover, as $w\leq
2\rho$,  we obtain from (\ref{last}) and (\ref{F(s,t)})
$$|\langle TP_{\mathcal E},F\Psi\rangle|\leq ||T||\cdot||\Psi||_{L_2}C^{1/2}\left(\sum_{i=1}^{n(\varepsilon)}
(\text{diam }\alpha_i)^{w}\right)^{1/2} \leq
(CK)^{1/2}||T||\cdot||\Psi||_{L_2}$$ so that
\begin{equation}\label{ineq1}
||F\cdot(TP_{\mathcal E})||_{{\mathfrak S}_2}\leq (CK)^{1/2}||T||
\end{equation}
Since $TP_{\mathcal E}\to^s T$, as $\e\to 0$, we have also
\begin{equation}\label{ineq2}
|\langle F\cdot T,\Psi\rangle|\leq D||T||||\Psi||_2
\end{equation}
 for any
$\Psi\in\Gamma(X,Y)\subset L_2(X\times Y,\mu\times\nu)$ and
$D=(CK)^{1/2}$. As $\Gamma(X,Y)$ is dense in $L_2(X\times Y,
\mu\times\nu)$, this implies $F\cdot T\in\mathfrak S_2$ with
\begin{equation}\label{ineq3}
||F\cdot T||_{\mathfrak S_2}\leq D||T||.
\end{equation}

 Now using  (\ref{ineq1}) and (\ref{ineq3}) we have
\begin{eqnarray*}
||(F\cdot T)P_{\mathcal E}-F\cdot(TP_{\mathcal E})||_{\mathfrak
S_2}&\leq& ||(F\cdot T)P_{\mathcal E}||_{\mathfrak S_2}
+||F\cdot(TP_{\mathcal
  E})||_{{\mathfrak S}_2}\\&\leq&
D||T||+D||T||=2D||T||.
\end{eqnarray*}
Since  $\langle (F\cdot T)P_{\mathcal E}-F\cdot(TP_{\mathcal
E}),\Psi\rangle\to 0$, as $\e\to 0$,  for any $\Psi\in\Gamma(X,Y)$,
we have that
 $(F\cdot T)P_{\mathcal E}-F\cdot(TP_{\mathcal E})\to 0$ weakly
in $\mathfrak S_2$.
 Then $$\displaystyle F\cdot T=\lim_{\varepsilon\to 0}(F\cdot T)P_{\mathcal E}=
\lim_{\varepsilon\to 0} ((F\cdot T)P_{\mathcal
E}-F\cdot(TP_{\mathcal E})+F\cdot(TP_{\mathcal E}))=
\lim_{\varepsilon\to 0} F\cdot(TP_{\mathcal
E})\in\overline{\text{Im}\Delta_F|_{\mathfrak
    S_2}}^w.$$

Thus $\Delta_F(T)\in\ker\Delta_F|_{\mathfrak
  S_2}\cap\overline{\text{Im}\Delta_F|_{\mathfrak S_2}}=\{0\}$.
\end{proof}

\section{Applications to linear operator equations in Hilbert space.}

Let ${\mathbb A}=\{A_k\}_{k\in K}$ be a commutative family of normal operators
with $\sum_{k\in K}||A_i||^2<\infty$. By $\sigma({\mathbb A})$ we denote the maximal
ideal space of the unital $C^*$-algebra generated by ${\mathbb A}$. To any
$t\in\sigma({\mathbb A})$ we associate  a sequence
$\lambda(t)=(t(A_1), t(A_2),\ldots)\in l_2$; the map
$t\mapsto\lambda(t)$ is continuous and identifies $\sigma({\mathbb
A})$ with a compact subset of $l_2$.

We say that the
{\it essential dimension}, $\text{ess-dim}$, of a ${\mathbb A} $ does
not exceed $r>0$ if there is a subset $D$ of $\sigma({\mathbb A})$
such that $E_{\mathbb A}(\sigma({\mathbb A})\setminus D)=0$ and
the Hausdorff dimension $\text{dim}(D)\leq r$.
In particular, if  all $A_k$ are Lipschitz
functions of one Hermitian (normal) operator with the Lipschitz
constants $C_k$ and  $\sum_{k\in K} C_k^2<\infty$ then
$\text{ess-dim}({\mathbb A})\leq 1$ (respectively $2$). If ${\mathbb
A}$ is diagonalizable then $\text{ess-dim}({\mathbb A})=0$.

Let for a bounded operator $\Delta: B(H_1,H_2)\to B(H_1,H_2)$,
${\mathcal E}_{\Delta}(0)$ be the root space of $\Delta$:
$${\mathcal E}_{\Delta}(0)=\{T\in B(H_1,H_2):
||\Delta^n(T)||^{1/n}\to 0, n\to\infty\}.$$

\begin{cor}\label{cor}
Let ${\mathbb A}=\{A_k\}_{k\in K}$, ${\mathbb B}=\{B_k\}_{k\in K}$
be two commutative families of normal operators on $H_1$ and $H_2$
satisfying $\sum_{k\in K}||A_k||^2<\infty$, $\sum_{k\in
K}||B_k||^2<\infty$. Suppose that we are given continuous functions
$f_j$, $g_j$ ( $j\in J$) on $\sigma({\mathbb A})$ and
$\sigma({\mathbb B})$ respectively and let $$\Delta(X)=\sum_{j\in
J}g_j({\mathbb B})Xf_j({\mathbb A})$$ for $X\in B(H_1,H_2).$ If
$\text{ess-dim } {\mathbb A}=w$ and $F(x,y)=\sum_{j\in
J}f_j(x)g_j(y)\in V_{\infty}(\sigma({\mathbb A}),\sigma(\mathbb B))$
satisfies (\ref{eq1}) with $E=\text{null }F$ and $\rho\geq w/2$ then
$\ker\Delta={\mathcal E}_{\Delta}(0)$.
\end{cor}

\begin{proof}
Assume  first that ${\mathbb A}$ and ${\mathbb B}$ have cyclic
vectors.  Then all $A_k$ and $B_k$ can be realized on
$L_2(\sigma({\mathbb A}),\mu)$ and $L_2(\sigma({\mathbb B}),\nu)$ as
multiplication operators by the coordinate functions and
$\Delta(X)=F\cdot X$, where $F(x,y)=\sum_{j\in J}f_j(x)g_j(y)$,
$(x,y)\in \sigma({\mathbb A})\times\sigma({\mathbb B})$. By
Theorem~\ref{opsyn} and \cite[Proposition~4.7]{sht2}, we obtain
$\ker\Delta_F={\mathfrak M}_{max}(\text{null }F)={\mathcal
E}_{\Delta_F}(0)$.

In the absence of a cyclic vector, decompose $H_1$ and $H_2$ into a direct sum
of subspaces
 $H_1=\oplus_{j=1}^{\infty}H_j^1$, $H_2=\oplus_{j=1}^{\infty} H_j^2$,
where each $H_j^1$ and $H_j^2$ is invariant with respect to
$\{A_k,A_k^*\}_{k\in K}$ and $\{B_k,B_k^*\}_{k\in K}$ respectively
and $\{A_k|_{H_j^1}\}_{k\in K}$, $\{B_k|_{H_j^2}\}$ has cyclic
vectors. Then each $X\in B(H_1, H_2)$ can be written as a
block-operator $X=(X_{ij})$, where $X_{ij}=P_{H_j^2}X|_{H_i^1}$ and
$P_{H_j^2}$ is the projection onto $H_j^2$, and
 $\Delta(X)=(\Delta_{ij}(X_{ij}))$,
where $\Delta_{ij}$ is the restriction of $\Delta$ to $B(H_i^1,H_j^2)$. Now
if $X\in{\mathcal E}_{\Delta}(0)$, $X_{ij}\in{\mathcal
E}_{\Delta_{ij}}(0)=\ker\Delta_{ij}$ and hence $X\in\ker\Delta$.
\end{proof}

\begin{cor}\label{power}
If, in notation of Corollary \ref{cor}, the function $F(x,y)$ satisfies (\ref{eq1}) with some $\rho > 0$ then
${\mathcal E}_{\Delta}(0) = \ker\Delta^k$ where $k = (w/2\rho]$, the least integer which is not less than $w/2\rho$.
\end{cor}
\begin{proof}
It suffices to note that the operator $\Delta^k$ corresponds to the function $F(x,y)^k$ and to apply Corollary \ref{cor}.
\end{proof}

\begin{cor}
Let ${\mathbb A}=\{A_k\}_{k\in K}$, ${\mathbb B}=\{B_k\}_{k\in K}$
be two commutative families of normal operators on $H_1$
and $H_2$ satisfying $\sum_{k\in K}||A_k||^2<\infty$, $\sum_{k\in
K}||B_k||^2<\infty$. Let $$\Delta(X)=\sum_{k\in
K}B_kXA_k \text{ for }X\in B(H_1,H_2).$$
If $\text{ess-dim } {\mathbb A}\leq 2n$  then
$\ker\Delta^n={\mathcal E}_{\Delta}(0)$.
\end{cor}
\begin{proof}
Follows from Corollary~\ref{power} with $f_k(x)=x_k$, $g_k(y)=y_k$,
for all $k\in K=J$.
\end{proof}

This result extends \cite[Theorem~10.3]{sht2} (where $n$ was actually equal 1) and simultaneously Corollary 10.4 of \cite{sht2} where the restriction $\text{card }K<\infty$ was present and essential for the proof.

\section{Applications to spectral synthesis in Varopoulos and Fourier algebras.}

 To formulate further results we have to recall the notion of spectral
synthesis in harmonic analysis.

Let ${\mathcal A}$ be a unital semi-simple regular commutative Banach algebra
with spectrum $X$, which is thus a compact Hausdorff space.
We will identify ${\mathcal A}$ with a subalgebra of the algebra $C(X)$ of
continuous
complex-valued functions on $X$ in our notation.
If $E\subseteq X$ is closed, let
\begin{align*}
I_{\mathcal A}(E)&=\{a\in{\mathcal A}:a(x)=0\text{ for }x\in E \}, \\
I_{\mathcal A}^0(E)&=\{a\in{\mathcal A}:a(x)=0\text{ in a nbhd of
}E\}
\end{align*}
 and
 $$J_{\mathcal A}(E)=\overline{I_{\mathcal A}^0(E)}.$$
One says that $a\in A$  {\it admits spectral synthesis} for ${\mathcal A}$ if
$a\in J_{\mathcal A}(\text{null }f)$.

Let $A^*$ be the dual of $A$.
 For $\tau\in A^*$ and
$a\in A$, define $a\tau$ in $A^*$ by $a\tau(b)=\tau(ab)$, and define
the support of $\tau$ by
$$\text{supp}(\tau)=\{x\in X_A: a\tau\ne 0\text{ whenever } a(x)\ne0\}.$$
In other words $\text{supp}(\tau)$ consists of all $x\in X_A$ such
that for any neighborhood $U$ of $x$ there exists $a\in A$ for which
$\text{supp}(a)\subset U$ and $\tau(a)\ne 0$. Then $a$ admits
spectral synthesis iff
$$\tau(a)=0 \text{ for each } \tau \text{ with }
\text{supp}(\tau)\subset\text{null}(a).$$

The Banach algebras we consider in this section are the Fourier
algebras $A({\mathbb T}^n )$ where ${\mathbb T}^n $ is the
$n$-dimensional torus, and the projective tensor product
$V(X,Y)=C(X)\hat\otimes C(Y)$, where $X$ and $Y$ are compact
Hausdorff spaces.

Note that $V(X,Y)$ (the Varopoulos algebra)
consists of all functions
$F\in C(X\times Y)$ which admit a representation
\begin{equation}\label{equa}
F(x,y)=\sum_{i=1}^{\infty}f_i(x)g_i(y),
\end{equation}
 where $f_i\in C(X)$, $g_i\in
C(Y)$ and $$\sum_{i=1}^{\infty}||f_i||_{C(X)}||g_i||_{C(Y)}<\infty.$$
$V(X,Y)$ is a Banach algebra  with respect to the pointwise multiplication and the norm
$$||F||_V=\inf\sum_{i=1}^{\infty}||f_i||_{C(X)}||g_i||_{C(Y)},$$
where $\inf$ is taken over all  representations of $F$ in the form
$\sum f_i(x)g_i(y)$ (shortly, $\sum f_i\otimes g_i$) satisfying the
above conditions (see \cite{var}).
It is known that $V(X,Y)$ is semi-simple, regular and its spectra is naturally identified with $X\times Y$.

As above we assume that $X$ is a metric space of Hausdorff dimension $w$.

\begin{theorem}\label{varopoulos}
Let $F\in V(X,Y)$  and a closed subset $E$ of $X\times Y$  satisfy
\begin{equation}\label{eq2}
|F(x,y)|\leq Cd(x,E^y)^{\rho}
\end{equation}
for all $(x,y)\in X\times Y$ and some $C>0$. If $\rho\geq w/2$ then
$F\in J_{V}(E)$. In particular, if (\ref{eq2}) holds for
$E=\text{null } F$ then $F$ admits spectral  synthesis in $V(X,Y)$.
\end{theorem}

\begin{proof}
By Theorem~\ref{opsyn}, for any choice of finite measures $\mu$ and
$\nu$,  $F\cdot T=0$ if  $T\in B(L_2(X,\mu),L_2(Y,\nu))$ such that
$\supp T\subset E$. The arguments in the proof of
\cite[Proposition~4.5]{sht2}  show that $F\in J_V(E)$.
\end{proof}

\begin{remark}\rm
For $w = n \in \mathbb N$, $X={\mathbb T}^{n}$, $Y={\mathbb T}^{m}$,
and $\rho > w/2$ the statement of Theorem~\ref{varopoulos}   was
proved by Varopoulos \cite[Theorem~7.2.2]{var} (using a different
method) under  a stronger condition on $F$:
\begin{equation}\label{eqvar}
|F(x,y)|\leq Cd((x,y),E)^{\rho}
\end{equation}
for all $(x,y)\in X\times Y$ and some $C>0$. Here $d$ is the natural
metric on ${\mathbb T}^{n}\times {\mathbb T}^{m}$. He also proved
that for any  $\rho < n/2$, one can construct a set $E\subset
{\mathbb T}^{n}\times {\mathbb T}^n$  and $F\in V({\mathbb T}^n,
{\mathbb T}^n)$ such that $F\in I_V(E)\setminus J_V(E)$; $|F(t)|\leq
C|t-E|^{\rho}$.

The case $\rho = n/2$ was left open in \cite{var}. It should be also
noted that the condition (\ref{eq2}) is less restrictive and more
convenient than (\ref{eqvar}), not only because it does not need any
metric on $Y$. For example, Corollary \ref{clasbp} below could not
be deduced from Theorem \ref{varopoulos} if in this theorem we had
(\ref{eqvar}) instead of (\ref{eq2}). This indicates that Theorem
\ref{varopoulos} is a natural form of the Beurling - Pollard theorem
for Varopoulos algebras.

\end{remark}
 The next result is the classical Beurling - Pollard Theorem for the Fourier algebras.
\begin{cor}\label{clasbp}
 If a function $F\in A({\mathbb T}^n)$ satisfies inequality
$$|f(t)|\leq Cd(t,\text{null } f)^{n/2}\quad\forall t\in {\mathbb T}^n, \text{
  some } C>0$$
then $f$ admits spectral synthesis in $A({\mathbb T}^n)$.
\end{cor}

\begin{proof}
Consider two linear mappings $M$ and $P$ (see \cite{var}):
$$A({\mathbb T}^n)\to^M V({\mathbb T}^n,{\mathbb T}^n)\to^P A({\mathbb T}^n)$$
defined for $f\in A({\mathbb T}^n)$ and $F\in V({\mathbb T}^n,
{\mathbb T}^n)$ by
$$Mf(x,y)=f(x+y), \quad PF(x)=\int_{{\mathbb T}^n}F(x-z,z)dz.$$
Then $M$ is an isometry and $P\circ M=Id_{A({\mathbb T}^n)}$.
For a closed set $E\subset {\mathbb T}^n$ let $E^*=\{(x,y): x+y\in
E\}$.
Then by \cite[Theorem~8.2.1]{var} $M^{-1}J_V(E^*)=J_A(E)$.
Thus to prove the statement it remains to show that $Mf\in
J_V((\text{null } f)^*)$.
This follows from the inequality
$$|Mf(x,y)|=|f(x+y)|\leq Cd(x+y,\text{null } f)^{\rho}=Cd(x,((\text{null } f)^*)^y)^{\rho}$$
and Theorem~\ref{varopoulos}.
\end{proof}

\section{Beurling-Pollard type theorem for weighted Fourier algebras.}\label{section5}

In what follows we need a modification of Hausdorff
dimension proposed in \cite{keckic}.
\begin{definition}\label{def}
For a compact metric space $X$, its {\it balanced Hausdorff dimension}, $bh(X)$, is
the infimum of all positive numbers $c$ that possess the following property:

 there exist positive constants $N$,
$P>0$ such that for all $\delta>0$ there exists a finite covering
$\displaystyle X\subset\sqcup_{j=1}^m\beta_j$
($\beta_i\cap\beta_j=\emptyset$) satisfying (i) $\delta/P<\text{diam
}\beta_j<\delta$ for all $1\leq j\leq m$ and (ii)
$\displaystyle\sum_{j=1}^m(\text{diam }\beta_j)^c\leq N$.
\end{definition}

Let $X$  be a compact metric space. For each ${\epsilon}>0$, let $N({\epsilon})$ be  the smallest number of balls of radius ${\epsilon}$ that cover $X$. Set
$$mo(X) = \liminf_{\epsilon\to 0}(- log N({\epsilon})/log {\epsilon}).$$
This characteristic ({\it the metric order} of $X$) was introduced
in \cite{PSh}.

\begin{proposition}\label{equality}
$mo(X) = bh(X)$.
\end{proposition}
\begin{proof}

Let $c > 0$. Clearly
\begin{equation}\label{ner-vo}
mo(X) < c \Longrightarrow N({\epsilon}) < Const~ {\epsilon}^{-c} \Longrightarrow mo(X)\le c.
\end{equation}

If $bh(X) < c$ then there are sets $\beta_1,...,\beta_n$ of
diameters between ${\epsilon}/P$ and ${\epsilon}$ such that they
cover $X$ and $\displaystyle\sum_{i=1}^n diam (\beta_i)^c < A$ (=
const). Hence $n \le A P^c/{\epsilon}^c$. Since any set of diameter
$ \le {\epsilon}$ is contained in a ball of radius $\le {\epsilon}$,
it follows that $N({\epsilon}) \le (AP^c)/ {\epsilon}^c$. So $mo(X)
\le c$.

Conversely if $mo(X) < c$, then let us take $N({\epsilon})$ balls
covering $X$ with radii $\le {\epsilon}$. Removing some parts of the
balls, we obtain non-intersecting sets $\beta_i$ of diameter $\le
{\epsilon}$ that cover $X$. We can then add to any of them a point
sufficiently far from it such that the diameter of each will be
${\epsilon}$. It will be a balanced covering (with $ P = 1$).
Now  $\sum_i diam(\beta_i)^c =
N({\epsilon}){\epsilon}^{c} < const$ by (\ref{ner-vo}).
\end{proof}

The result shows, in particular, that if $f:X\to Y$ is a surjective Lipschitz map between metric spaces then $bh(Y)\le bh(X)$.

Another easy consequence of Proposition \ref{equality} is the next lemma, proved in \cite[Lemma~2.2]{keckic} by other method.

Let us denote, for any $\varepsilon > 0$ and any subset $E$ of a
metric space $X$, by $E_{\varepsilon}$ the
$\varepsilon$-neighborhood of $E$, i.e., $E_{\varepsilon}=\{x\in X:
d(x,E)\leq \varepsilon\}$.
\begin{lemma}\label{lemmakeckic}
Let $E\subset{\mathbb R}^{n}$ be a set of balanced Hausdorff
dimension $<c$ and let for all $\e>0$ $E_{\e}$ be the open
$\e$-neighborhood of $E$. Then $m(E_{\e})\leq D\e^{n-c}$ for some
constant $D>0$.
\end{lemma}
\begin{proof}
Given $\e>0$, let $\displaystyle E=\bigcup_{i=1}^{N(\e)}B_{i,\e}$ be
a covering of $E$ by the smallest number, $N(\varepsilon)$, of balls
of radius $\e$. Clearly, $\displaystyle
E_{\e}\subset\bigcup_{i=1}^{N(\e)}B_{i,2\e}$, where $B_{i,2\e}$ is
the ball  with the same centrum as $B_{i,\e}$ and the radius $2\e$.
Let $S_{n}$ denote the unit ball in ${\mathbb R}^{n}$. Then, as
$bh(E)=mo(E)<c$, $N(\e)\leq \text{const}\e^{-c}$, giving
$$m(E_{\e})\leq\sum_{i=1}^{N(\e)}m(B_{i,2\e})=m(S_{n})N(\e)(2\e)^{n}\leq
\text{const}\e^{n-c}.$$
\end{proof}

Let us for any $\alpha>0$, denote by $A_{\alpha}({\mathbb T}^n)$ the
algebra of functions $f\in C({\mathbb T}^n)$ such that
$$f(t)=\sum_{k\in{\mathbb Z}^n}a_ke^{ik\cdot t} \text{ with }
||f||_{A_{\alpha}}=\sum_{k\in{\mathbb
Z}^n}|a_k|(1+|k|)^{\alpha}<\infty,$$ where for $k=(k_1,\ldots,
k_n)$, $t=(t_1,\ldots,t_n)$ we write $k\cdot t=k_1t_1+\ldots k_nt_n$
and $|k|=(k_1^2+\ldots +k_n^2)^{1/2}$. Let $PM_{\alpha}({\mathbb
T}^n)$ denote the dual space of $A_{\alpha}({\mathbb T}^n)$.

\begin{theorem}\label{balan}
Let $E\subset{\mathbb T}^n$ be a closed set of balanced Hausdorff
dimension $c$, and let $f\in
I_{A_{\alpha}}(E)$ be such that $|f(x)|\leq Ad(x,E)^m$ for any $x\in
{\mathbb T}^n$ and some $A>0$, $m> 0$. If $m>c/2+\alpha$ or
$m=c/2+\alpha=n/2+\alpha$, then $f\in J_{A_{\alpha}}(E)$.
\end{theorem}

\begin{proof}
For any $\varepsilon>0$ and $\beta > 0$, set
$$\delta(x)=\left\{\begin{array}{cc}
\displaystyle(1-|x|^2)^{\beta},&|x|\leq 1,\\
\displaystyle 0,& |x|>1
\end{array}\right. $$
and $\delta_{\e}(x)=\e^{-n}\delta(x/\e)/||\delta_1||_1$
($||\cdot||_1$ is the $L_1({\mathbb T}^n)$-norm). Then
$||\delta_{\e}||_1=1$ and we have
$$\hat\delta_{\varepsilon}(k)=\int_{{\mathbb
    T}^n}\delta_{\e}(t)e^{ikt}\to
e^{ik\cdot 0}=1, \text{ as }\e\to 0.$$

Thus, for $f\in A_{\alpha}({\mathbb T}^n)$, $T\in
PM_{\alpha}({\mathbb T}^n)$,
$$\langle T,f\rangle=\lim_{\e\to 0}\sum_{k\in{\mathbb Z}^n}\hat
T(k)\hat f(k)\hat\delta_{\e}(k)=\lim_{\e\to 0}
(T*\delta_{\e},f),$$
where $(\cdot,\cdot)$ is the scalar product in $L_2({\mathbb T}^n)$.

Our next step is to estimate  $||T*\delta_{\e}||_2$:
\begin{eqnarray*}
||T*\delta_{\e}||_{L_2}&=&||\{\hat
T(k)\hat\delta_{\e}(k)\}||_{l_2}\leq ||\{\hat
T(k)/(1+|k|)^{\alpha}\}||_{l_{\infty}}||\{\hat\delta_{\e}(k)(1+|k|)^{\alpha}\}||_{l_2}\\&=&
||T||_{PM_{\alpha}}(\sum_{k:|k|\leq\frac{1}{\e}}|\hat\delta_{\e}(k)|^2(1+|k|)^{2\alpha}+
\sum_{k:|k|>\frac{1}{\e}}|\hat\delta_{\e}(k)|^2(1+|k|)^{2\alpha})^{1/2}.
\end{eqnarray*}

It is easy to see that $\text{sup}_{\e,k}|\hat\delta_{\e}(k)|=C<\infty$. Hence
$$\sum_{k:|k|\leq\frac{1}{\e}}|\hat\delta_{\e}(k)|^2(1+|k|)^{2\alpha}\leq
C\sum_{k:|k|\leq\frac{1}{\e}}(1+|k|)^{2\alpha}\leq
D\sum_{|k|=0}^{\left[\frac{1}{\e}\right]}|k|^{2\alpha}|k|^{n-1}\leq
D\frac{1}{\e^{2\alpha+n-1}}\cdot\frac{1}{\e}=\frac{D}{\e^{2\alpha+n}}$$
for some other constant $D$.

To estimate the second sum we use the explicit formulas for the Fourier
coefficients $\delta_{\e}(k)$, given in
\cite[Theorem~4.15]{stein_weiss}:
$\hat\delta_{\e}(k)=C(\beta)|k\e|^{-n/2-\beta}J_{n/2+\beta}(2\pi|k\e|)$,
where $J_{\nu}$ is the Bessel function and $C(\beta)$ is a constant
depending on $\beta$.

As $|J_{\nu}(r)|\leq C_{\nu}r^{-1/2}$ for $r\geq 1$,
this gives us that for $\e|k|>1$
$$|\hat\delta_{\e}(k)|\leq C(\beta)|k\e|^{-n/2-1/2-\beta}$$
for some other constant $C(\beta)$. Thus, choosing $\beta>\alpha-1/2$,
we obtain from the following estimation that $\delta_{\e}\in
A_{\alpha}({\mathbb T}^n)$ and
\begin{eqnarray*}
\sum_{k:|k|>\frac{1}{\e}}|\hat\delta_{\e}(k)|^2(1+|k|)^{2\alpha}&\leq&
\frac{C(\beta)^2}{\e^{n+1+2\beta}}\sum_{k:|k|>\frac{1}{\e}}|k|^{2\alpha-n-1-2\beta}\\&\leq&
\frac{C(\beta)^2}{\e^{n+1+2\beta}}\sum_{|k|=\left[\frac{1}{\e}\right]+1}^{\infty}|k|^{2\alpha-n-1-2\beta}|k|^{n-1}\\&\leq&
\frac{C(\beta)^2}{\e^{n+1+2\beta}}\int_{\left[\frac{1}{\e}\right]}^{\infty}x^{2\alpha-2-2\beta}
dx\\&\leq&
\frac{D}{\e^{n+1+2\beta}}\cdot\frac{1}{\e^{2\alpha-1-2\beta}}=\frac{D}{\e^{2\alpha+n}}.
\end{eqnarray*}
Let $\supp T\subset E$. Then $\supp T*\delta_{\e}\subset E_{\e}$ and
\begin{eqnarray*}
|\langle T*\delta_{\e},f\rangle|&\leq&
||T*\delta_{\e}||_{L_2}\left(\int_{E_{\e}}|f(x)|^2dx\right)^{1/2}\\&\leq&
||T||_{PM_{\alpha}}||\{\hat\delta_{\e}(k)(1+|k|)^{\alpha}\}||_{l_2}\sup_{E_{\e}}|f|(m(E_{\e}\setminus\text{null
  }f))^{1/2}\\&\leq& \frac{D}{\e^{n/2+\alpha}}||T||_{PM_{\alpha}}\sup_{E_{\e}}|f|(m(E_{\e}\setminus\text{null
  }f))^{1/2}
\end{eqnarray*}
for some constant $D$.

 If $m>\alpha+c/2$, we choose $\gamma>0$ so
that $m>\alpha+c/2+\gamma$. As $E$ has balanced Hausdorff dimension
$c<c+2\gamma$ then, by Lemma~\ref{lemmakeckic},
$m(E_{\e}\setminus\text{null }f)\leq m(E_{\e})\leq
C\e^{n-c-2\gamma}$.
 As $|f(x)|\leq
A~d(E,x)^m$ for any $x$, $y\in {\mathbb T}^n$, we obtain now
$$|\langle T*\delta_{\e}, f\rangle|\leq
D\frac{1}{\e^{n/2+\alpha}}\e^m\e^{n/2-c/2-\gamma}\leq
D\e^{m-c/2-\alpha-\gamma}.$$ Thus  we obtain
$$\langle T,f\rangle=\lim_{\e\to 0}\langle
T*\delta_{\e},f\rangle=0.$$

If $m=c/2+\alpha=n/2+\alpha $ we have also
$$|\langle T*\delta_{\e},f\rangle|\leq D\frac{1}{\e^{n/2+\alpha}}\e^m(m(E_{\e}\setminus\text{null
  }f))^{1/2}=D(m(E_{\e}\setminus\text{null
  }f))^{1/2}\to 0,\text{ as }\e\to 0,$$
and $\langle T,f\rangle=0$.
\end{proof}

\begin{cor}\label{beurling}
Let $E\subset{\mathbb T}^n$ be a closed subset of balanced Hausdorff
dimension $c$, and let $f\in I_A(E)$ be such that $|f(x)|\leq
Ad(x,E)^m$ for any  $x\in {\mathbb T}^n$ and some $A>0$,  $m>0$. If
$m>c/2$ or $m=n/2=c/2$ then $f\in J_A(E)$.
\end{cor}

\section{Elementary operators and spectral synthesis}

Now we apply the results of Section~\ref{section5} to the
study of ascent of some operators on Banach spaces and, in particular, elementary operators on Banach algebras.

Let $B$ be a semisimple, regular, commutative Banach algebra with
unit and let $M$ be a Banach $B$-module. For any $x\in M$ set
\begin{eqnarray*}
&ann(x)=\{b\in B\mid b\cdot x=0\},\\
&Supp(x)= null(ann(x)).
\end{eqnarray*}
Then $ann(x)$ is a closed ideal in $B$ and $Supp(x)$ is a closed
subset in $X_B$, the spectrum of $B$. More generally, for any subset
${\mathcal L} \subset M$ we denote by $Supp({\mathcal L})$ the
smallest closed set $Supp({\mathcal L})$ such that $Supp(x)\subset
Supp({\mathcal L})$ for any $x\in{\mathcal L}$. We will need the
following simple lemma.

\begin{lemma}\label{lemma}
Let $B$ be a commutative semisimple regular Banach algebra and let
$M$ be a Banach $B$-module. If $f\in J_B(E)$ then  $f\cdot x=0$ for
any $x\in M$ such that $Supp (x)\subset E$.
\end{lemma}

\begin{proof}
In fact, $f\in J_B(E)\subset J_B( Supp (x))\subset ann(x)$.
\end{proof}

Let $A$ be a unital (not necessarily commutative) Banach algebra. We
say that an element $a\in A$ is a {\it  generalized scalar element}
if there are constants $C$ and $s$ such that \
\begin{equation}\label{gen}
||e^{ita}||\leq C(1+|t|)^s,\forall t\in{\mathbb R}.
\end{equation}
In this case we say that $a$ has order $s$.

In particular, we may speak about {\it generalized scalar
operators}. The concept has its origin in the work of Colojoara and
Foias (see \cite{cf}).

Let $T = (T_1,...,T_n)$ be a commutative family of generalized scalar operators on a Banach space $\mathcal X$.
We call by the order of $T$ the number ${\alpha} = s_1 +...+s_n$ where $s_i$ are the orders of $T_i$.

We denote by $\sigma(T)$ the joint spectrum of the family $T$ (it
was shown in \cite{mprr} that all classical definitions of joint
spectrum coincide for commutative families of generalized scalar
operators). We denote by $c$ the balanced Hausdorff dimension of
$\sigma(T)$.

The family $T$ admits an $A_{\alpha}({\mathbb T}^n)$-calculus:
$$f(T) = \sum_{m\in {\mathbb Z}^n}\hat{f}(m)\exp(imT)$$
where $mT = \sum m_iT_i$ for $m = (m_1,...,m_n)$.

For a linear operator $\Delta$ on ${\mathcal X}$ we will mean by the
{\it ascent}, $\text{asc }\Delta$, of $\Delta$ the least positive
integer such that $\ker\Delta^m=\ker\Delta^{m+1}$.

\begin{theorem}\label{genasc}
If $g\in A_{\alpha}({\mathbb T}^n)$ is a Lipschitz function, and $V
= g(T)$ then $\text{asc }V \le [c/2+\alpha]+1$.
\end{theorem}
\begin{proof}
We introduce an $A_{\alpha}({\mathbb T}^n)$-module structure on
$\mathcal X$ setting $f\cdot x = f(T)x$, $x\in {\mathcal X}$. It is
not difficult to see that $Supp({\mathcal X}) = \tau(T):=
\{(\exp(t_1),...,\exp(t_n)): (t_1,...,t_n)\in \sigma(T)\}$. Hence
$Supp(x) \subset \tau(T)$ for any $x\in \mathcal X$.

Let $V^Nx = 0$ for some $N$.  Then $g^N\in ann(x)$, $Supp(x)\subset
null(g^N) = null(g)$. Thus setting $E = \tau(T)\cap null(g)$ we get
that $Supp(x)\subset E$.

Since $bh(E)\le c$ and $|g^m(t)| \le C~d(t,E)^m$ for each $m$, it
follows from Theorem \ref{balan} that $g^m\in J(E)$   if $m > c/2 +
\alpha$. Applying Lemma \ref{lemma} we see that $g^m\cdot x = 0$.
This means that $V^mx = 0$.
\end{proof}

Let now $a_1$, $a_2,\ldots, a_n$ and $b_1$, $b_2,\ldots,b_n$ be
commuting $n$-tuples of generalized scalar elements of a unital
Banach algebra $A$ of orders $s_1,\ldots, s_n$ and $r_1,\ldots, r_n$
Let $s=s_1+\ldots s_n$ and $r=r_1+\ldots r_n$. For $a\in A$ we
denote by $L_a$ and $R_a$ the operators acting on $A$ via the left
and, respectively, the right multiplication by $a$. Let
$T=(L_{a_1},\ldots L_{a_n}, R_{b_1},\ldots R_{b_n})$; $T$ is a
commuting family of generalized scalar operators and the order  of
$T$ does not exceed $\alpha = s+r$.

We consider the "elementary" operator
$\displaystyle\Lambda=\sum_{i=1}^nL_{a_i}R_{b_i}$ on $A$.

\begin{cor}\cite{keckic}
Assume that the balanced Hausdorff dimension of $\sigma(T)$ is less
than or equal to $c$. Then $asc(\Lambda)\leq
[c/2+\alpha]+1$.
\end{cor}
\begin{proof}
We may assume that the norms of all $a_i, b_i$ are less than $K <
\pi$. Then $\Lambda = g(T)$ where
$g(t)=\sum_{i=1}^{n}\varphi(t_i)\varphi(t_{n+i})$ and $\varphi(t)$
is a smooth $2\pi$-periodic function such that $\varphi(t)=t$ for
$|t|< K$. Hence our statement follows from Theorem \ref{genasc}.
\end{proof}

\vspace{0.2cm}

\begin{center}
{\bf Acknowledgment}
\end{center}
The work was partially written when the first author was visiting
Chalmers University of Technology in G\"oteborg, Sweden. The
research  was partially supported by a grant from the Swedish Royal
Academy of Sciences as a part of the program of cooperation with the
former Soviet Union. The second author was also supported by the
Swedish Research Council.

\begin{center}
{\footnotesize \sl Department of Mathematics, Vologda State Liceum of
Mathematical and Natural Sciences,
Vologda, 160000, Russia}\\
{\footnotesize \sl
 shulman\_v\symbol{64}yahoo.com}\\
 \vspace{0.5cm}
{\footnotesize \sl Department of Mathematics, Chalmers University of
Technology,
SE-412 96 G\"oteborg, Sweden} \\
{\footnotesize \sl
turowska\symbol{64}math.chalmers.se}
\end{center}

\end{document}